\documentclass[a4paper,11pt]{article}

\usepackage{MyLaTeX,MyArticle}
\input xy
\xyoption{all}

\title{Control of Fusion and Solubility in Fusion Systems}
\author{David A.~Craven}
\date{February 2009}
\renewcommand{\F}{\mathcal F}
\newcommand{\E}{\mathcal E}
\renewcommand{\centre}[1]{\mathrm{Z}(#1)}

\begin{document}
\maketitle

\begin{abstract} In this article, we consider control of fusion, quotients, and $p$-soluble fusion systems. For control of fusion, we prove the three main theorems in the literature in a new, largely elementary way, significantly shortening their proofs. To prove one of these, and a theorem of Aschbacher that the product of strongly closed subgroups is strongly closed, we produce a consolidated treatment of quotients, collating and expanding the constructions previously available. We move on to $p$-soluble fusion systems, and prove that they are constrained, allowing us to effectively characterize fusion systems of $p$-soluble groups. This leads us to recast Thompson Factorization for $Qd(p)$-free fusion systems, and consider it for more general fusion systems.
%In this article, we consider the control of fusion in fusion systems, proving three previously known, non-trivial results in a new, largely elementary way. We then reprove a result of Aschbacher, that the product of two strongly closed subgroups is strongly closed; to do this, we consolidate the theory of quotients of fusion systems into a consistent theory. We move on considering $p$-soluble fusion systems, and prove that they are constrained, allowing us to effectively characterize fusion systems of $p$-soluble groups. This leads us to recast Thompson Factorization for $Qd(p)$-free fusion systems, and consider Thompson Factorization for more general fusion systems.
\end{abstract}

\section{Introduction}

The theory of fusion systems, which was started by Puig but lay unpublished for decades, attempts to formalize the concept of local finite group theory, and extend such results to the sphere of blocks of finite groups. It has attracted considerable and growing attention over the course of the last decade, but several basic questions remain unresolved in the field. In addition, some of the results in the literature permit considerable shortenings in their proofs; the improved exposition will make clear the reasons behind the results.

\begin{thma}[Stancu {{\cite[Proposition 6.2]{stancu2006}}}]\label{normaliffnormal} Let $\F$ be a saturated fusion system on a finite $p$-group $P$, and let $Q$ be a subgroup of $P$. Then $\F=\Norm_\F(Q)$ if and only if $\F_Q(Q)\normal\F$.
\end{thma}

The proof given here is the first proof that does not rely on a non-trivial result: in \cite{stancu2006}, this result was given as a corollary to what is here Theorem \ref{normaliffstrongcentral}; in \cite{linckelmann2007}, this result was a corollary of what is here Theorem \ref{normalcontrol}; and in \cite{aschbacher2007un}, the proof of this result requires the non-trivial fact that a constrained fusion system is the fusion system of a finite group (see \cite[Proposition C]{bcglo2005}). Here we will prove it by elementary means, using the extension of the Frattini argument to fusion systems, as proved in \cite{aschbacher2008}.

Using Theorem \ref{normaliffnormal}, we derive Theorem \ref{normalcontrol} as a consequence.

\begin{thma}[Linckelmann {{\cite[Theorem 9.1]{linckelmann2007}}}]\label{normalcontrol} Let $\F$ be a saturated fusion system on a finite $p$-group $P$, and let $\E$ be a normal subsystem, on a subgroup $Q$. Suppose that $\E=\Norm_\E(R)$ for some subgroup $R$, and let $S$ be the subgroup of $Q$ generated by all $\F$-conjugates of $R$. Then $\F=\Norm_\F(S)$.
\end{thma}

Finally, using an easy result of Aschbacher's on equivalent conditions to fusion being controlled by a subgroup, we derive Stancu's main theorem of \cite{stancu2006}.

\begin{thma}[Stancu {{\cite[Theorem 4.8]{stancu2006}}}] \label{normaliffstrongcentral} Let $\F$ be a saturated fusion system on a finite $p$-group $P$, and let $Q$ be a subgroup of $P$. Then $\F=\Norm_\F(Q)$ if and only if $Q$ is strongly $\F$-closed, and there is a central series
\[1=Q_0\leq Q_1\leq \cdots \leq Q_n=Q\]
with each of the $Q_i$ weakly $\F$-closed.
\end{thma}

%Moving away from control of fusion -- i.e., a subgroup $Q$ such that $\F=\Norm_\F(Q)$ -- we will reprove two theorems of Aschbacher. The first has a highly non-trivial proof in \cite{aschbacher2007un}.

Moving away from control of fusion -- i.e., a subgroup $Q$ such that $\F=\Norm_\F(Q)$ -- we will reprove a theorem of Aschbacher. It has a highly non-trivial proof in \cite{aschbacher2007un}.

\begin{thma}[Aschbacher {{\cite[Theorem 2]{aschbacher2007un}}}]\label{productstronglyclosed} Let $\F$ be a fusion system. Then the product of two strongly $\F$-closed subgroups is strongly $\F$-closed.
\end{thma}

The proof of this result in \cite{aschbacher2007un} is difficult, and requires considerable preliminaries. Here we will produce an extremely short and trivial proof of the result. However, in order to do so, we must consider factor systems, a subject that has not yet been fully understood, with differing accounts of it in the literature (see \cite{puig2006}, \cite{stancu2006}, \cite{linckelmann2007}, and \cite{aschbacher2008} for four, all different to various extents, approaches). Here we combine the various approaches to produce the first comprehensive treatment of the subject in Section \ref{quotients}.

We end with a discussion of $p$-soluble systems, which has not been considered in the literature before, although there are some unpublished notes on the subject by various authors. Broadly speaking, a fusion system is called \emph{$p$-soluble} if repeated quotienting out by $\Orth_p(\F)$ eventually exhausts the group: a formal definition will be given in Section \ref{psoluble}. The main result in this section is the following. (In the context of the generalized Fitting subsystem, this is also proved in \cite{aschbacher2007un}.)

\begin{thma} Let $\F$ be a $p$-soluble fusion system; then
\[ \Cent_P(\Orth_p(\F))\leq \Orth_p(\F).\]
In particular, $\F$ is a constrained fusion system, so is the fusion system of a finite group.
\end{thma}

In light of the fact that saturated subsystems of $p$-soluble systems are $p$-soluble, this gives an easy proof of the known fact that all block fusion systems of $p$-soluble groups are group fusion systems.

In the final section, we consider the fusion systems of $p$-soluble groups; these can be effectively characterized internally.

\begin{thma}\label{intcharsol} Let $\F$ be a saturated fusion system, and write $Q=\Orth_p(\F)$. Then $\F$ is the fusion system of a $p$-soluble group if and only if $\F$ is constrained and $\Aut_\F(Q)$ is $p$-soluble.
\end{thma}

Using this, we make some remarks on extending Thompson Factorization to constrained fusion systems, and how such results differ from the versions for groups.

We begin in the next section with the definitions and concepts from fusion systems that we need. Because the notation and terminology has not yet crystallized, we will define everything to avoid ambiguity; this section contains Alperin's fusion theorem. After a preliminary section we set about proving the first two theorems given above. In the final sections, from Section \ref{quotients} onwards, we introduce and study factor systems, prove Theorems \ref{normaliffstrongcentral} and \ref{productstronglyclosed}, and then consider $p$-soluble systems. Note that our maps will be composed from left to right.

\section{Fusion Systems}

Here we collect the very basic concepts in the field of fusion systems. We begin with defining fusion systems, and before that, a prefusion system, which is needed to define fusion systems before one proves that they are, in fact, fusion systems.

\begin{defn} Let $P$ be a finite $p$-group. Then a \emph{prefusion system} $\F$ consists of the set of all subgroups of $P$, and for each pair of subgroups $Q$ and $R$, a set $\Hom_\F(Q,R)$ of injective homomorphisms $Q\to R$. The composition of two morphisms is given by the composition of group homomorphisms, when it lies in the prefusion system.

A \emph{fusion system} $\F$ on $P$ is a prefusion system, with composition of morphisms $\phi:Q\to R$ and $\psi:R\to S$ always possible in $\F$, and whose morphisms $\Hom_\F(Q,R)$ should satisfy the following three axioms:
\begin{enumerate}
\item for each $g\in P$ with $Q^g\leq R$, the associated conjugation map $\theta_g:Q\to R$ is in $\Hom_\F(Q,R)$;
\item for each $\phi\in\Hom_\F(Q,R)$, the isomorphism $Q\to Q\phi$ lies in $\Hom_\F(Q,Q\phi)$; and
\item if $\phi\in\Hom_\F(Q,R)$ is an isomorphism, then its inverse $\phi^{-1}:R\to Q$ lies in $\Hom_\F(R,Q)$.
\end{enumerate}
We write $\Aut_\F(Q)$ for the set (in fact group) $\Hom_\F(Q,Q)$.
\end{defn}

If $G$ is a finite group with Sylow $p$-subgroup $P$, we write $\F_P(G)$ for the fusion system got from conjugation by elements of $G$ on the subgroups of $P$. We also note the \emph{universal fusion system} $\mc U(P)$ on a finite $p$-group $P$, where $\Hom_{\mc U(P)}(Q,R)$ is the set of \emph{all} injective homomorphisms from $Q$ to $R$. Obviously any fusion system is a subsystem of the universal fusion system.

\begin{defn} Let $P$ be a finite $p$-group, and let $Q$ be a subgroup of $P$. Let $\F$ be a fusion system on $P$. We say that $Q$ is \emph{fully normalized} if, whenever $\phi:Q\to R$ is an isomorphism in $\F$, we have that $|\Norm_P(Q)|\geq |\Norm_P(R)|$.
\end{defn}

\begin{defn} Let $P$ be a finite $p$-group, and let $\F$ be a fusion system on $P$. We say that $\F$ is \emph{saturated} if
\begin{enumerate}
\item $\Aut_P(P)$ is a Sylow $p$-subgroup of $\Aut_\F(P)$, and
\item every morphism $\phi:Q\to P$ in $\F$ such that $Q\phi$ is fully normalized extends to a morphism $\bar\phi: N_\phi\to P$, where
\[ N_\phi=\{x\in \Norm_P(Q)\,:\,\text{there exists }y\in \Norm_P(Q\phi)\text{ such that }(g^x)\phi=(g\phi)^y\text{ for all }g\in Q\}.\]
\end{enumerate}
\end{defn}

Note that $N_\phi$ is the inverse image under the map $\delta_Q:\Norm_P(Q)\to\Aut(Q)$ of the subgroup $\Aut_P(Q)\cap \phi\Aut_P(R)\phi^{-1}$, and that $Q\Cent_P(Q)$ is always contained within $N_\phi$.

\begin{defn} Let $\F$ be a fusion system on the finite $p$-group $P$. Let $Q$ be a subgroup of $P$, and let $K$ be a subgroup of $\Aut(Q)$.
\begin{enumerate}
\item The subgroup $\Norm_P^K(Q)$, the \emph{$K$-normalizer} of $Q$ in $P$, is the set of all $g\in \Norm_P(Q)$ such that $g$ induces an automorphism of $Q$ that lies in $K$. Write $\Aut_P(Q)$ for the subgroup of $\Aut(Q)$ of automorphisms induced by elements of $P$, and $\Aut_P^K(Q)=K\cap\Aut_P(Q)$.
\item
The fusion system $\Norm_\F^K(Q)$ is the category whose objects are all subgroups of $\Norm_P^K(Q)$, and whose morphisms $\Hom_{\Norm_\F^K(Q)}(R,S)$ are
\[ \{\phi\in\Hom_\F(R,S)\,:\,\phi\text{ extends to }\bar\phi\in\Hom_\F(QR,QS)\text{ with }\bar\phi|_Q\in K\}.\]
The fusion system $\Norm_\F^K(Q)$ is called the \emph{$K$-normalizer subsystem} in $\F$ of $Q$.
\end{enumerate}
Of particular interest are the cases $K=\Aut(Q)$, in which case we write $\Norm_\F(Q)$ and call this the \emph{normalizer}, and $K=1$, in which case we write $\Cent_\F(Q)$ and call this the \emph{centralizer}.
\end{defn}

If $Q$ is fully normalized and $K\normal\Aut_\F(Q)$ then $\Norm_\F^K(Q)$ is saturated (see \cite[Theorem 3.6]{linckelmann2007}). If $Q$ is a subgroup of $P$, we are interested in the case where $\F=\Norm_\F(Q)$, in which case we say that the subgroup is \emph{normal}. Theorems \ref{normaliffnormal} and \ref{normaliffstrongcentral} give equivalent conditions for a subgroup to be normal. We need a few more of these in the course of this article, but in order to state the first one in Section \ref{prelims} we need the concept of centric and radical subgroups.

Let $\F$ be a fusion system on a finite $p$-group $P$, and let $Q$ be a subgroup of $P$. We say that $Q$ is \emph{$\F$-centric} if, whenever $R$ is $\F$-isomorphic to $Q$, then $R$ contains its centralizer (or equivalently, $\Cent_P(R)=\centre R$). We say that $Q$ is \emph{$\F$-radical} if $\Orth_p(\Aut_\F(Q))=\Inn(Q)$. A system is \emph{constrained} if it contains a normal, $\F$-centric subgroup. The fundamental theorem on such systems is the following.

\begin{thm}[{{\cite[Proposition C]{bcglo2005}}}]\label{constrainedmodels} Let $\F$ be a constrained fusion system on a finite $p$-group $P$. Then there is a unique finite group $G$ such that 
\begin{enumerate}
\item $O_{p'}(G)=1$,
\item $\Cent_G(\Orth_p(G))\leq \Orth_p(G)$, and
\item $\F=\F_P(G)$.
\end{enumerate}
\end{thm}

Using the concept of centric and radical subgroups, we also have a version for fusion systems of Alperin's fusion theorem.

\begin{thm}[Alperin's fusion theorem] Let $\F$ be a saturated fusion system on a finite $p$-group $P$, and let $\phi: Q\to R$ be an isomorphism. Then there exist
\begin{enumerate}
\item a sequence of $\F$-isomorphic subgroups $Q=Q_0$, $Q_1$, \dots, $Q_n=R$,
\item a sequence $S_1$, $S_2$, \dots, $S_n$ of fully normalized, $\F$-radical, $\F$-centric subgroups, with $Q_{i-1},\; Q_i\leq S_i$, and
\item a sequence of $\F$-automorphisms $\phi_i$ of $S_i$ such that $Q_{i-1}\phi_i=Q_i$,
\end{enumerate}
such that
\[ (\phi_1\phi_2\ldots \phi_n)|_Q=\phi.\]
\end{thm}

A proof may be found in, for example, \cite{linckelmann2007}.

Now we define two types of subgroups: a subgroup $Q$ is called \emph{weakly $\F$-closed} if $Q$ is only $\F$-isomorphic to itself, and it is called \emph{strongly $\F$-closed} if whenever any $\F$-morphism has domain a subgroup of $Q$, it has image a subgroup of $Q$.

\begin{defn}Let $\F$ be a fusion system on a finite $p$-group $P$, and let $\E$ be a subsystem on a subgroup $Q$ of $P$, where $Q$ is strongly $\F$-closed. We say that $\E$ is \emph{$\F$-invariant} if, for each $R\leq S\leq Q$, $\phi\in\Hom_{\E}(R,S)$, and $\psi\in \Hom_{\F}(S,P)$, we have that $\psi^{-1}\phi\psi$ is a morphism in $\Hom_{\E}(R\psi,Q)$. If, in addition, $\E$ is saturated, we say that $\E$ is \emph{normal} in $\F$. We denote normality by $\E\normal\F$. A saturated fusion system is \emph{simple} if it contains no non-trivial, proper, normal subsystems.

If $\F$ is a saturated fusion system and $\E$ is a normal subsystem, then $\E$ is said to be \emph{characteristic} if $\E\phi=\E$ for all $\phi\in\Aut(\F)$. We denote this by $\E\charac \F$.
\end{defn}

Because we need to use a variant of the Frattini argument for fusion systems, we need the definition of an $\F$-Frattini subsystem.

\begin{defn}[{{\cite[Section 3]{aschbacher2008}}}] Let $\F$ be a fusion system on a finite $p$-group $P$, and let $\E$ be a subsystem of $\F$, on the subgroup $Q$. We say that $\E$ is \emph{$\F$-Frattini} if, whenever $R\leq Q$ and $\phi:R\to P$ is a morphism in $\F$, there exist morphisms $\alpha\in\Aut_\F(Q)$ and $\beta\in\Hom_\E(R\alpha,P)$ such that $\phi=\alpha\beta$.
\end{defn}

\section{Preliminaries}

\label{prelims}

Our first result tells us that we may apply the Frattini argument whenever $\E$ is a normal subsystem of $\F$, and in fact gives us an equivalent condition to normality of subsystems.

\begin{thm}[Aschbacher {{\cite[Theorem 3.3]{aschbacher2008}}}]\label{normalifffrattini} Let $\F$ be a saturated fusion system on a finite $p$-group $P$, and let $\E$ be a subsystem on a strongly $\F$-closed subgroup $Q$. Then $\E$ is $\F$-invariant if and only if $\Aut_\F(Q)\leq \Aut(\E)$ and $\E$ is $\F$-Frattini.
\end{thm}

We also need the equivalent for fusion systems of the statement that if $K\charac H\normal G$ then $K\normal G$.

\begin{prop}[Aschbacher {{\cite[7.4]{aschbacher2007un}}}]\label{charnormal} Let $\F$ be a saturated fusion system on a finite $p$-group $P$. Suppose that $\E\normal \E'\normal \F$, and write $Q$ for the subgroup on which $\E'$ acts. Suppose that all $\F$-automorphisms of $Q$ induce automorphisms on $\E$. Then $\E\normal \F$. In particular, if $\E\charac \E'\normal \F$, then $\E\normal \F$.
\end{prop}
\begin{pf} Suppose that $\E\charac \E'$. Then $\Aut_\F(Q)\leq \Aut(\E')$ by Theorem \ref{normalifffrattini}, and so since $\E$ is $\Aut(\E')$-invariant, the second statement follows from the first. Let $R$ be the subgroup of $Q$ on which $\E$ acts.

If $S\leq R$ and $\phi:S\to P$ is any map in $\F$, then since $Q$ is strongly $\F$-closed, we have that $\im\phi\leq Q$. Since $\E'\normal \F$, it is $\F$-Frattini by Theorem \ref{normalifffrattini}, and so $\phi=\alpha\beta$, where $\alpha\in\Aut_\F(Q)$ and $\beta$ is a morphism in $\E'$. Since all $\F$-automorphisms of $Q$ induce automorphisms on $\E$, we see that $S\alpha\leq R$, and as $R$ is strongly $\E'$-closed, $S\phi=(S\alpha)\beta\leq R$. In particular, $R$ is strongly $\F$-closed as $S\phi\leq R$.

We will use Theorem \ref{normalifffrattini} again: since $\E\normal \E'$, $\E$ is $\E'$-Frattini, so $\beta=\gamma\delta$, where $\gamma\in\Aut_{\E'}(R)$ and $\delta$ is a morphism in $\E'$. Then $\phi=\alpha\gamma\delta$. We claim that this decomposition proves that $\E$ is $\F$-Frattini; since $\alpha\in\Aut_\F(Q)$ and $R$ is strongly $\F$-closed, $\alpha|_R$ is an automorphism, and hence $\alpha\gamma\in\Aut_\F(R)$. Since $\delta\in\Hom_\E(S\alpha\gamma,R)$, this gives the correct decomposition of $\phi$. 

If $S=R$, then we perform the same decomposition to get that $\phi\in\Aut_\F(R)$ may be written as $\phi=\alpha\gamma$, where $\alpha\in\Aut_\F(Q)$ and $\gamma\in\Aut_{\E'}(R)$. Because all elements of $\Aut_\F(Q)$ induce automorphisms on $\E$ by hypothesis, and $\Aut_{\E'}(R)\leq \Aut(\E)$ by Theorem \ref{normalifffrattini}, we see that $\Aut_\F(R)\leq \Aut(\E)$, and so $\E\normal \F$ by Theorem \ref{normalifffrattini}.
\end{pf}

We also need two equivalent conditions for a subgroup to be normal in a fusion system.

\begin{prop}[{{\cite[Proposition 1.6]{bcglo2005}}}]\label{normalifffrc} Let $\F$ be a saturated fusion system on a finite $p$-group $P$, and let $Q$ be a strongly $\F$-closed subgroup of $P$. Then $\F=\Norm_\F(Q)$ if and only if $Q$ is contained in every fully normalized, $\F$-centric, $\F$-radical subgroup.
\end{prop}
%\begin{pf} By \cite[Proposition 5.6]{linckelmann2007}, if $\F=\Norm_\F(Q)$, then $Q$ is contained in every $\F$-radical, $\F$-centric subgroup. Therefore, suppose that $Q$ is contained in every fully normalized, $\F$-radical, $\F$-centric subgroup, and let $\phi:R\to S$ be a morphism in $\F$. By Alperin's fusion theorem, $\phi$ maybe the decomposed as the restriction of $\phi_1\phi_2\ldots\phi_n$, where $\phi_i\in\Aut_\F(U_i)$ for fully normalized, $\F$-radical, $\F$-centric subgroups $U_i$. Let $R=R_0$ and $R_i=R_{i-1}\phi_i$, so that $R_{i-1}$ is $\F$-conjugate to $R_i$ via $\phi_i$ for all $1\leq i\leq n$. Since $Q$ is contained in each $U_i$, we have restrictions $\psi_i:QR_{i-1}\to QR_i$ of each $\phi_i$, and the composition of the $\psi_i$ gives $\psi:QR\to QS$, extending $\phi$. Since $Q$ is strongly $\F$-closed, $\psi|_Q\in\Aut_\F(Q)$, as required.
%\end{pf}

\begin{prop}[Aschbacher]\label{normaliffstrongclosedcentral} Let $\F$ be a saturated fusion system on a finite $p$-group $P$, and let $Q$ be a subgroup of $P$. Then the following are equivalent:
\begin{enumerate}
\item $\F=\Norm_\F(Q)$; and
\item there exists a central series
\[ 1=Q_0\leq Q_1\leq \cdots \leq Q_n=Q\]
for $Q$, all of whose terms are strongly $\F$-closed.
\end{enumerate}
\end{prop}
\begin{pf} Suppose that $\F=\Norm_\F(Q)$, and let $Q_i=\central iQ$. We need to show that every $Q_i$ is strongly $\F$-closed, and we are done. Let $R\leq Q_i$, and $\phi\in\Hom_\F(R,P)$. Then $\phi$ extends to $\bar\phi\in\Hom_\F(Q,P)$ that acts like an automorphism on $Q$, and since $Q_i\charac Q$, we see that $\im \phi\leq Q_i$, as needed.

Now suppose that (ii) holds, and let $T$ be any $\F$-radical, $\F$-centric subgroup. If we can show that in this case $Q\leq T$, then we are done, since then $Q$ is contained in every $\F$-radical, $\F$-centric subgroup, and so $\F=\Norm_\F(Q)$ by Proposition \ref{normalifffrc}. Choose $i$ maximal such that $Q_i\leq T$, so that $Q_{i+1}\not\leq T$. (For a contradiction, assume that $Q$ is not a subgroup of $T$.) Set $R=Q_{i+1}\cap T$ and $S=\Norm_{Q_{i+1}}(T)$; as $Q\not\leq T$, we have that $R<S$. Since $S$ normalizes $T$ and $Q_{i+1}\normal P$, we have that $[T,S]\leq R$, and as $Q_{i+1}/Q_i$ is central in $Q/Q_i$, we see that $S$ centralizes $R/Q_i$, and since each $Q_{j+1}/Q_j$ is central in $Q/Q_j$, we have that $S$ centralizes $Q_{j+1}/Q_j$ for all $j<i$.

Each $Q_i$ is strongly $\F$-closed, and therefore $\Aut_\F(T)$ acts on $R=Q_{i+1}\cap T$ and $Q_j$ for all $j\leq i$. Hence there is an $\Aut_\F(T)$-invariant  series
\[ 1=Q_0\leq \cdots \leq Q_i\leq R\leq T,\]
with $S$ centralizing each factor. The set of all such automorphisms of $T$ is clearly a normal subgroup, and is a $p$-subgroup by \cite[Corollary 5.3.3]{gorenstein}. Therefore, $\Aut_S(T)$ is contained in a normal $p$-subgroup of $\Aut_\F(T)$, and so in particular $\Aut_S(T)\leq \Inn(T)$ since $T$ is $\F$-radical. Thus $S=\Norm_{Q_{i+1}}(T)\leq T\Cent_P(T)=T$ since $T$ is $\F$-centric. Therefore, $Q_{i+1}\cap T\geq S>R=Q_{i+1}\cap T$, a contradiction. Hence $Q$ is contained in every $\F$-centric, $\F$-radical subgroup, and so $\F=\Norm_\F(Q)$, as claimed.
\end{pf}

We are now in a position to prove three preliminary lemmas -- the last one perhaps of independent interest -- needed to prove the first three main theorems.

\begin{lem}[cf.\ {{\cite[Proposition 6.2]{stancu2006}}}]\label{controlimpnormal} Let $\F$ be a saturated fusion system, and let $Q$ be a subgroup such that $\F=\Norm_\F(Q)$. Then $\F_Q(Q)$ is normal in $\F$.
\end{lem}
\begin{pf} Suppose that $\F=\Norm_\F(Q)$; then any morphism $\phi:R\to S$ lifts to a morphism $\bar\phi:QR\to QS$ that acts as an automorphism on $Q$ and so $Q$ is strongly $\F$-closed. Write $\E=\F_Q(Q)$; if $R\leq S\leq Q$ and $\psi:S\to P$ is a map in $\F$, then we should show that, for all $\phi=\theta_g$ with $g\in Q$, the map $\psi^{-1}\phi\psi$ is also in $\E$. If $\psi:S\to P$ is a map, then since $Q$ is strongly $\F$-closed, this is actually $\psi:S\to Q$, and since $\F=\Norm_\F(Q)$, the map $\psi$ extends to an automorphism $\bar\psi\in\Aut_\F(Q)$. Thus we may assume that $S=Q$. If $\phi=\theta_g\in\Hom_\E(R,Q)$ is a morphism for some $g\in Q$, then it clearly extends to a map $\theta_g\in\Aut_\E(Q)$, and so we may assume that $R=Q$ as well. Thus suppose that $R=S=Q$; then we need to show that if $g\in Q$ and $\psi\in\Aut_\F(Q)$, then $\psi^{-1}\theta_g\psi\in\F_Q(Q)$, but it is well-known that $\psi^{-1}\theta_g\psi=\theta_{g\psi}\in \F_Q(Q)$, and we get the result.
\end{pf}

The proof of the next lemma is due to Sejong Park, and replaces a less elegant proof in a preliminary version of the manuscript.

\begin{lem}\label{fqqnormalimpstrongcentral} Let $\F$ be a saturated fusion system, and let $Q$ be a subgroup such that $\F_Q(Q)\normal \F$. Then every characteristic subgroup of $Q$ is strongly $\F$-closed. In particular, $Q$ has a central series each of whose terms is strongly $\F$-closed.
\end{lem}
\begin{pf} Let $R$ be a characteristic subgroup of $Q$. Then $\F_R(R)\charac \F_Q(Q)\normal \F$, and so $\F_R(R)\normal \F$ by Proposition \ref{charnormal}; in particular, $R$ is strongly $\F$-closed.
\end{pf}

%\begin{pf} Let $R$ be a characteristic subgroup of $Q$, and note that since $\E=\F_Q(Q)\normal \F$, we must have that $Q$ is strongly $\F$-closed. Let $S$ be any subgroup of $R$, and let $\phi:S\to T$ be an $\F$-isomorphism. By Theorem \ref{normalifffrattini}, $\phi$ may be factorized as $\phi=\alpha\beta$, where $\alpha\in\Aut_\F(Q)$ and $\beta\in\Hom_\E(S\alpha,T)$. Since $R$ is characteristic, and by Theorem \ref{normalifffrattini} $\Aut_\F(Q)\leq \Aut(\E)$, we have that $S\alpha\leq R\alpha=R$. Finally, since characteristic subgroups are clearly strongly $\F_Q(Q)$-closed, we must have that 
%\[ S\phi=(S\alpha)\beta\leq R,\]
%as claimed.
%\end{pf}

\begin{lem}\label{weakcentralimpstrong} Let $\F$ be a saturated fusion system on a finite $p$-group $P$, and let $Q$ be a strongly $\F$-closed subgroup of $P$. Let $Z$ be a central subgroup of $Q$ that is weakly $\F$-closed. Then $Z$ is strongly $\F$-closed.
\end{lem}
\begin{pf} Let $\phi:X\to Y$ be an $\F$-isomorphism with either $X$ or $Y$ contained in $Z$, and $Y$ fully normalized. Since $Q$ is strongly $\F$-closed, both $X$ and $Y$ lie inside $Q$. The subgroup $Z$ is central in $Q$, and so $Z\leq \Cent_P(X)$; hence $\phi$ extends to a morphism $\psi:ZX\to Q$ which, as $Z$ is weakly $\F$-closed, restricts to an automorphism $\psi|_Z:Z\to Z$. Therefore if either $X$ or $Y$ lies in $Z$, then the other also does, as required.
\end{pf}
% Let $Y$ be a subgroup of $Z$, and let $\phi:Y\to R$ be an $\F$-isomorphism such that $R$ is fully normalized. Since $\F$ is saturated, and $Z\leq Y\Cent_P(Y)\leq N_\phi$ always, we have that $\phi$ extends to a map $\psi:Z\to Q$. Since $Z$ is weakly $\F$-closed, the image of this map must also be $Z$, and so $R\leq Z$. Hence all fully normalized subgroups $\F$-conjugate to $Y$ lie inside $Z$, and so we assume that $Y$ is fully normalized itself. Let $\chi:S\to Y$ be any $\F$-isomorphism, and note that since $Q$ is strongly $\F$-closed, we have that $S\leq Q$. Again, $Z\leq N_\chi$ since $Z$ is central in $Q$, and so $\chi$ extends to $\theta:SZ\to Q$. Since $Z$ is weakly $\F$-closed, this means that $Z\leq \im \theta$ as well. Since $Z$ contains $Y$, we must have that $S=Y\theta^{-1}\leq Z\theta^{-1}=Z$, and so $Z$ is strongly $\F$-closed.

\section{Theorems \ref{normaliffnormal} and \ref{normalcontrol}}

We begin by proving Theorem \ref{normaliffnormal}.

\begin{thm} Let $\F$ be a saturated fusion system on a finite $p$-group $P$, and let $Q$ be a subgroup of $P$. Then $\F=\Norm_\F(Q)$ if and only if $\F_Q(Q)\normal\F$.
\end{thm}
\begin{pf} This follows immediately from Proposition \ref{normaliffstrongclosedcentral} and Lemmas \ref{controlimpnormal} and \ref{fqqnormalimpstrongcentral}.\end{pf}

We now prove Theorem \ref{normalcontrol}, using Theorem \ref{normaliffnormal}.

\begin{thm} Let $\F$ be a saturated fusion system on a finite $p$-group $P$, and let $\E$ be a normal subsystem, on a subgroup $Q$. Suppose that $\E=\Norm_\E(R)$ for some subgroup $R$, and let $S$ be the subgroup of $Q$ generated by all $\F$-conjugates of $R$. Then $\F=\Norm_\F(S)$.
\end{thm}
\begin{pf} Suppose that $\E=\Norm_\E(R)$, and write $S$ for the product of all $\F$-conjugates of $R$ (all of which are subgroups of $Q$, since $Q$ is strongly $\F$-closed). Then $\E=\Norm_\E(S)$ since the product of subgroups normal in $\E$ is also normal in $\E$, and if $\phi\in\Aut_\F(Q)$, then $\phi$ leaves $S$ invariant, so induces an automorphism on $\mc D=\F_S(S)$. (By Theorem \ref{normaliffnormal}, $\F_S(S)\normal \F$ if and only if $\F=\Norm_\F(S)$.) Therefore, since $\mc D\normal \E\normal \F$, by Proposition \ref{charnormal}, $\mc D\normal\F$. Again, by Theorem \ref{normaliffnormal}, we get that $\F=\Norm_\F(S)$, as claimed.
\end{pf}

Related to this, we have a proposition on the subgroup $\Orth_p(\F)$ of a fusion system $\F$. As we said earlier, if $Q$ and $R$ are normal subgroups of $\F$, then $QR$ is as well, and hence there is a largest subgroup of $P$ normal in $\F$; this is denoted by $\Orth_p(\F)$. Denote by $\mc O_p(\F)$ the subsystem $\F_Q(Q)$ of $\F$, where $Q=\Orth_p(\F)$. By Theorem \ref{normaliffnormal}, $\mc O_p(\F)$ is a normal subsystem, and it is invariant under $\Aut(\F)$, so characteristic.

\begin{prop}\label{normalop} Let $\F$ be a saturated fusion system on a finite $p$-group $P$. If $\E$ is a normal subsystem of $\F$ on a subgroup $Q$, then $\Orth_p(\F)\cap Q=\Orth_p(\E)$. In particular, if $\Orth_p(\E)\neq 1$ for some $\E\normal\F$, then $\Orth_p(\F)\neq 1$.
\end{prop}
\begin{pf} Let $R=\Orth_p(\F)$; by Proposition \ref{normaliffstrongclosedcentral}, $R$ possesses a central series
\[ 1=R_0\leq R_1\leq \cdots \leq R_d=R,\]
such that each $R_i$ is strongly $\F$-closed. We claim that $Q_i=Q\cap R_i$ is strongly $\E$-closed; in this case,
\[ 1=Q_0\leq Q_1\leq \cdots \leq Q_d=Q\cap R\]
is a central series for $Q\cap R$, each of whose terms is strongly $\E$-closed, yielding that $Q\cap R\leq \Orth_p(\E)$, as required. It remains to show that $Q_i$ is strongly $\E$-closed; however, any morphism in $\E$ that originates inside $Q_i=Q\cap R_i$ must have image inside $Q$ since $\E$ lies on $Q$, and must also lie in $R_i$ since it is strongly $\F$-closed, and so $Q_i$ is strongly $\E$-closed. Thus $Q\cap R\leq \Orth_p(\E)$.

On the other hand, Theorem \ref{normalcontrol} tells us that $\Orth_p(\E)\leq \Orth_p(\F)$, and so we get equality.
\end{pf}

In order to prove Theorems \ref{normaliffstrongcentral} and \ref{productstronglyclosed}, we need to understand factor systems, and we delay their proofs until after the next section.

\section{Quotients}
\label{quotients}

In this section we will consider morphisms of fusion systems and quotients. The treatment of these varies wildly in the literature with several opposing viewpoints and one or two errors, and it is our intention here to produce a clear description of the subject. We begin by defining a factor system.

\begin{defn} Let $Q$ be a normal subgroup of $P$, and let $\F$ be a fusion system on $P$. By the \emph{factor system $\F/Q$}, we mean the prefusion system on $P/Q$, such that for any two subgroups $R$ and $S$ containing $Q$, we have that $\Hom_{\F/Q}(R/Q,S/Q)$ is the set of homomorphisms $\phi$ induced from the set $\Hom_\F(R,S)$ such that $Q\phi=Q$.
\end{defn}

Traditionally, in the definition above the subgroup $Q$ is strongly $\F$-closed, but this is not necessary for the following two results. The first is easy, and a proof is omitted.

\begin{prop} Let $Q$ be a normal subgroup of $P$, and let $\F$ be a fusion system on $Q$. The prefusion system $\F/Q$ is a fusion system on $P/Q$.
\end{prop}
%\begin{pf} If $\phi:R/Q\to S/Q$ and $\psi:S/Q\to T/Q$ are morphisms in $\F/Q$, then there are morphisms $\phi':R\to S$ and $\psi':S\to T$ in $\F$ that are preimages of $\phi$ and $\psi$ respectively, and $\phi'\psi':R\to T$ has image $\phi\psi$ in $\F/Q$. Hence $\F/Q$ is a category.

%Certainly $\F_{P/Q}(P/Q)$ is contained within $\F/Q$, since conjugation by a coset $Qx$ on $P/Q$ is the same as that induced by $x$ on $P/Q$. If $\phi:R/Q\to S/Q$ is a morphism in $\F/Q$, then it is induced by a morphism $\phi':R\to S$ with $Q\phi'=Q$. As $\F$ is a fusion system, the corresponding $\F$-isomorphism $\psi':R\to R\phi$ lies in $\F$, and since $Q\psi'=Q$, $Q$ lies inside both $R$ and $R\psi'$. The second axiom of a fusion system is satisfied by $\F/Q$ because $\psi:R/Q\to (R/Q)\psi$ is induced by $\psi':R\to R\psi'$. Finally, if $\phi':R\to S$ and $\psi':S\to R$ are $\F$-isomorphisms with $Q\leq R,S$, $Q\phi'=Q\psi'=Q$, and $\phi'\psi'=1$, then the induced morphisms $\phi$ and $\psi$ are mutually inverse as well, proving that $\F/Q$ is, indeed, a fusion system.
%\end{pf}

Our proof of the next proposition follows \cite[Theorem 6.2]{linckelmann2007}, and although our hypotheses are weaker, the method of proof is the same.

\begin{prop}\label{quotientissaturated} Let $\F$ be a saturated fusion system on a finite $p$-group $P$, and let $Q$ be a weakly $\F$-closed subgroup of $P$. Then the fusion system $\F/Q$ is saturated.
\end{prop}
\begin{pf} All automorphisms in $\Aut_{\F/Q}(P/Q)$ are induced from automorphisms in $\Aut_\F(P)$, and so the obvious homomorphism $\Aut_\F(P)\to \Aut_{\F/Q}(P/Q)$ is surjective. The image of $\Aut_P(P)$ in $\Aut_{\F/Q}(P/Q)$ is clearly $\Aut_{P/Q}(P/Q)$, so that it satisfies the first axiom of a saturated fusion system.

Suppose that $\phi\in\Hom_{\F/Q}(R/Q,S/Q)$ is an isomorphism such that $S/Q$ is fully $\F/Q$-normalized. We claim that $S$ is also fully $\F$-normalized. Since $Q\leq R$, and $Q$ is weakly $\F$-closed, for all $T$ that are $\F$-isomorphic to $R$, we have that $Q\leq T$ and $Q\leq \Norm_P(T)$. Also, $\Norm_P(R)/Q=\Norm_{P/Q}(R/Q)$. Therefore
\[ |\Norm_P(T)|=|\Norm_{P/Q}(T/Q)|\cdot|Q|\leq |\Norm_{P/Q}(S/Q)|\cdot|Q|=|\Norm_P(S)|;\]
hence $S$ is fully $\F$-normalized.

\bigskip

Now let $\phi$ be an automorphism of a fully $\F/Q$-normalized subgroup $R/Q$, and let $\psi$ be an $\F$-automorphism of $R$ with image $\phi$ in $\F/Q$. At this point we would like to prove that $N_\psi/Q=N_\phi$, but it is not necessarily true. However, there is \emph{some} $\psi$ with image $\phi$ for which $N_\psi/Q=N_\phi$, as we shall demonstrate now. Notice that $N_\psi/Q\leq N_\phi$ trivially.

Let $K$ be the kernel of the natural map $\Aut_\F(R)\to \Aut(R/Q)$, a normal subgroup of $\Aut_\F(R)$; then $K$ consists of all elements of $\Aut_\F(R)$ that act trivially on $R/Q$, and hence are sent to the identity automorphism of $R/Q$ under the map $\F\to\F/Q$. The idea is that if $\chi\in K$, then $\chi\psi$ and $\psi$ both have the image $\phi$ in $\F/Q$, so one may `ignore' elements in $K$. We will prove that there are morphisms $\chi\in K$ and $\theta:R\to R$ such that $\psi=\chi\theta$ and $\theta$ has the property that $N_\theta/Q=N_\phi$. Since $\theta$ and $\psi$ define the same image $\phi$ in $\F/Q$, we prove that $\phi$ extends to $N_\phi$.

Since $K$ is a normal subgroup of $\Aut_\F(R)$ and $\Aut_P(R)$ is a Sylow $p$-subgroup of $\Aut_\F(R)$ (as $R$ is fully normalized), we have that $\Aut_P^K(R)=K\cap \Aut_P(R)$ is a Sylow $p$-subgroup of $K$, and by the Frattini argument
\[ \Aut_\F(R)=K\Norm_{\Aut_\F(R)}(\Aut_P^K(R)).\]
Since $\Aut_P(R)$ normalizes $\Aut_P^K(R)=\Aut_P(R)\cap K$, we see that $\Aut_P(R)$ is a Sylow $p$-subgroup of $\Norm_{\Aut_\F(R)}(\Aut_P^K(R))$. Also, writing $X=\Norm_{\Aut_\F(R)}(\Aut_P^K(R))$, we have
\[ X/X\cap K\cong KX/K=\Aut_\F(R)/K\cong \Aut_{\F/Q}(R/Q),\]
by the second isomorphism theorem and the definition of $K$. Since $X=\Norm_{\Aut_\F(R)}(\Aut_P^K(R))$, we may form the quotient group $X/\Aut_P^K(R)$, and as $S=\Aut_P^K(R)$ is a Sylow $p$-subgroup of $K$, it must be a (normal) Sylow $p$-subgroup of $X\cap K$; hence $(X/\cap K)/S$ is a $p'$-group. Quotienting out by this Sylow $p$-subgroup, we see that
\[ X/X\cap K\cong (X/S)/(X\cap K)/S \cong \Aut_{\F/Q}(R/Q).\]

Notice that $N_\phi$ is the preimage in $P/Q$ of the intersection $\Aut_{P/Q}(R/Q)\cap \Aut_{P/Q}(R/Q)^{\phi^{-1}}$ of two Sylow $p$-subgroups of $\Aut_{\F/Q}(R/Q)$. Since $(X\cap K)/S$ is a $p'$-group, we see that there are two Sylow $p$-subgroups $A/S$ and $B/S$ of $X/S$ that project onto $\Aut_{P/Q}(R/Q)$ and $\Aut_{P/Q}(R/Q)^{\phi^{-1}}$ respectively, and an element $Sg\in X/S$, such that $A/S\cap (B/S)^{Sg^{-1}}$ projects onto the correct intersection and $Sg$ projects onto $\phi$. (For a proof, see \cite[Corollary 2.13]{linckelmann2007}; it essentially comes from the fact that if $H$ is a normal $p'$-subgroup, then the fusion systems of $G$ and $G/H$ are isomorphic.)

Since $S$ is a normal $p$-subgroup, we may pull this up to $X$, and so there is some element $\theta\in X=\Norm_{\Aut_\F(R)}(\Aut_P^K(R))$ such that $\Aut_P(R)\cap \Aut_P(R)^{\theta^{-1}}$ maps onto $\Aut_{P/Q}(R/Q)\cap (\Aut_{P/Q}(R/Q))^{\phi^{-1}}$. Therefore the map $N_\theta\to N_\phi$ is surjective, and so $N_\theta/Q=N_\phi$. As the map $\theta$ extends to $N_\theta$, the map $\phi$ extends to $N_\phi$.

\bigskip

It remains to deal with any map $\phi:S/Q\to R/Q$ in $\F/Q$, where $R/Q$ is fully $\F/Q$-normalized. This lifts to a map $\psi:S\to R$ in $\F$ with $R$ fully $\F$-normalized. By \cite[Lemma 2.6]{linckelmann2007}, there is some map $\theta:S\to R$ with $N_\theta=N_P(S)$, and so $\phi$ extends to $N_\phi$ if and only if both the image $\chi$ of $\theta$ in $\F/Q$ extends to $N_\chi=\Norm_{P/Q}(S/Q)$ and $\chi^{-1}\phi$ extends to $N_{\chi^{-1}\phi}=N_\phi$. The first of these claims is obvious, and the second has been proved, concluding the proof.
\end{pf}

In the other direction, since fusion systems are categories, one may consider morphisms of fusion systems.

\begin{defn} Let $\F$ and $\E$ be fusion systems on the finite $p$-groups $P$ and $Q$ respectively. Then a \emph{morphism} $\Phi:\F\to\E$ of fusion systems  is a pair $(\phi,\{\phi_{R,S}:R,S\leq P\})$, where $\phi:P\to Q$ is a group homomorphism, and for each $R,S\leq P$, the map $\phi_{R,S}$ is a function
\[ \phi_{R,S}:\Hom_\F(R,S)\to \Hom_\E(R\phi,S\phi)\]
such that the corresponding map $\F\to\E$ on the category forms a functor; i.e., for any two composable $\F$-morphisms $\alpha$ and $\beta$, we have
\[ (\alpha\beta)\Phi=\alpha\Phi\beta\Phi.\]
\end{defn}

It is easy to see, using the functoriality of $\Phi$, that the action of $\Phi$ in morphisms in $\F$ is completely determined by the underlying group homomorphism.

Kernels of morphisms are clearly normal subgroups, since they are kernels of group homomorphisms. We have even more.

\begin{prop} Let $\F$ be a fusion system on the finite $p$-group $P$. Let $\phi$ be a morphism of fusion systems from $\F$. Then $\ker \phi$ is strongly $\F$-closed.
\end{prop}
\begin{pf} Let $Q$ be the kernel of $\phi$, and let $R$ be a subgroup of $Q$. We need to show that if $S$ is $\F$-isomorphic to $R$ then $S\leq Q$. Let $\psi:R\to S$ be an isomorphism. Then $\psi\phi$ is an isomorphism in the target fusion system, and since $S\phi$ is trivial, we must have that $R\phi$ is trivial also. Thus $Q$ is strongly $\F$-closed.
\end{pf}

We will now construct, for every strongly $\F$-closed subgroup $Q$, a morphism on $\F$ with kernel exactly $Q$.

\begin{defn} Let $P$ be a finite $p$-group and let $\F$ be a fusion system on $P$. Let $Q$ be a strongly $\F$-closed subgroup of $P$. By $\bar\F_Q$, we will denote the prefusion system on $P/Q$ with morphisms $\Hom_{\bar\F_Q}(R/Q,S/Q)$ consisting of those morphisms induced by $\Hom_\F(R',S')$, as $R'$ and $S'$ range over all subgroups of $P$ such that $R'Q=R$ and $S'Q=S$. (Since $Q$ is strongly $\F$-closed, any such morphism $\phi:R'\to S'$ gives rise to a morphism $\bar\phi:R'Q/Q\to S'Q/Q$.)

By $\gen{\bar\F_Q}$ we denote the prefusion system on $P/Q$ consisting of all finite compositions of morphisms from $\bar\F_Q$.
\end{defn}

Notice that $\F/Q$ is contained inside $\bar\F_Q$. It turns out that $\gen{\bar\F_Q}$ is a fusion system.

\begin{lem} Let $Q$ be a strongly $\F$-closed subgroup of a finite $p$-group $P$, on which a fusion system $\F$ is defined. Then $\gen{\bar\F_Q}$ is a fusion system on $P/Q$.
\end{lem}
\begin{pf} That $\gen{\bar\F_Q}$ is a category is obvious, since we are guaranteed compositions of morphisms by definition. The first axiom of a fusion system is satisfied, since
\[ \F_{P/Q}(P/Q)\subs \F/Q\subs \gen{\bar\F_Q}.\]
If we prove the final two axioms for the subset $\bar\F_Q$, then since $\gen{\bar\F_Q}$ is got from $\bar\F_Q$ by compositions of morphisms, those axioms would also hold there. If $\phi:R/Q\to S/Q$ is a morphism in $\bar\F_Q$, then there is some morphism $\phi':R'\to S'$ inducing $\phi$, and the corresponding isomorphism $\psi':R'\to R'\phi'$ induces an isomorphism in $\bar\F_Q$ corresponding to $\phi$. Finally, if $\phi:R/Q\to S/Q$ is an isomorphism in $\bar\F_Q$, then it comes from some isomorphism $\phi':R'\to S'$ in $\F$, and the inverse of $\phi'$ induces the inverse of $\phi$. Hence $\bar\F_Q$ is a fusion system.
\end{pf}

At this point it becomes difficult to know whether to define the natural map $\F\to \gen{\bar\F_Q}$ as a morphism of fusion systems, since although it satisfies the requirements in the definition, the image of the map, $\bar\F_Q$, is not in general a fusion system, because it is not a category. If $\bar\F_Q=\gen{\bar\F_Q}$, then the natural map does become a surjective morphism of fusion systems.

\begin{prop} Let $\F$ be a fusion system on a finite $p$-group $P$, and suppose that $Q$ is a subgroup such that $\F=\Norm_\F(Q)$. Then $\F/Q=\bar\F_Q$, and hence the natural map $\F\to\F/Q$ is a morphism of fusion systems.
\end{prop}
\begin{pf} Any morphism $\phi:R\to S$ extends to a morphism $\psi:QR\to QS$ that acts as an automorphism on $Q$. Certainly, $\bar\phi=\bar\psi$ in $\bar\F_Q$, since the action of $\phi$ and $\psi$ on $QR/Q$ is the same. Also, $\bar\psi\in\F/Q$, and since $\F/Q\subs \bar\F_Q$ we must have equality.
\end{pf}

In general of course, a strongly $\F$-closed subgroup need not be normal in $\F$, and in this case we need not have that $\bar\F_Q$ is a fusion system, or that $\gen{\bar\F_Q}=\F/Q$.

\begin{example} Let $P=\gen{a,b,c,d}$ be elementary abelian of order $16$, and let $\F$ be the fusion system generated by $\F_P(P)$ and the two morphisms $\gen{ab}\to\gen{c}$ and $\gen{ac}\to \gen{d}$. Then $A=\gen a$ is strongly $\F$-closed, and so we may form the prefusion system $\bar\F_A$. Here, the cosets $Ab$ and $Ac$ are $\bar\F_A$-conjugate, as are the cosets $Ac$ and $Ad$. However, there is no map sending $Ab$ to $Ad$, and so $\bar\F_A$ is not a fusion system.

Note also that there are no non-trivial morphisms on overgroups of $A$, and so $\F/A=\F_{P/A}(P/A)$.
\end{example}

Other than the case where $\F=\Norm_\F(Q)$, there is another case in which $\F/Q=\bar\F_Q$. The proof of this theorem follows \cite[Proposition 6.3]{puig2006}, although simplifications have been made.

\begin{thm}\label{saturatedquotient} Let $\F$ be a saturated fusion system and let $Q$ be a strongly $\F$-closed subgroup. Then $\F/Q=\bar\F_Q$, and so the map $\F\to\F/Q$ is a morphism of fusion systems.
\end{thm}
\begin{pf} If $\phi:R\to S$ is a map in $\F$, write $\bar\phi$ for the image of this map in $\bar\F_Q$; i.e., write $\bar\phi$ for the induced map $\bar\phi:QR/Q\to QS/Q$. Firstly, notice that both $\F/Q$ and $\bar\F_Q$ are on the same subgroup, namely $P/Q$. Certainly, $\F/Q$ is contained in $\bar\F_Q$, so we need to prove the converse; in other words, given a morphism $\phi:R\to S$ in $\F$, we need to show that there is some morphism $\psi:RQ\to P$ such that $\bar\phi=\bar\psi$, for then $\bar\phi\in\F/Q$, as needed.

We proceed by induction on $n=|P:R|$, noting that if $R$ contains $Q$ then we are done trivially; in particular, this implies that $n>1$. By Alperin's fusion theorem, any morphism may be factored as (restrictions of) a sequence of automorphisms $\phi_i$ of fully normalized, $\F$-centric, $\F$-radical subgroups $U_i$. Suppose that the images $\bar\phi_i$ of each of the $\phi_i$ lie in $\F/Q$: since $\F/Q$ is a fusion system and $\F\to\gen{\bar\F_Q}$ is a morphism, we have that $\bar\phi$, the product of (restrictions of) the $\bar\phi_i$, also lies in $\F/Q$.

Therefore we may assume that one of the $\bar\phi_i$ lies in $\bar\F_Q$ but not in $\F/Q$. By our inductive hypothesis, we see that $|U_i|=|R|$, and so we may replace $R$ and $\phi$ by $U_i$ and $\phi_i$; therefore $R$ is now a fully normalized, $\F$-centric, $\F$-radical subgroup, and $\phi$ is an automorphism of $R$. Also, $Q\not\leq R$ as we saw above.

Similar to the proof of Proposition \ref{quotientissaturated}, let $K$ be the kernel of the natural map $\Aut_\F(R)\to \Aut(QR/R)$, a normal subgroup of $A=\Aut_\F(R)$; then $K$ consists of all elements of $\Aut_\F(R)$ that act trivially on $QR/Q$, and hence are sent to the identity automorphism of $RQ/Q$ under the map $\F\to\F/Q$.

Let $T=\Norm_P^K(R)$; since $R$ is fully normalized, $\Aut_P(R)$ is a Sylow $p$-subgroup of $A$, and so $K\cap \Aut_P(R)=\Aut_T(R)$ is a Sylow $p$-subgroup of $K$. Therefore, by the Frattini argument,
\[ A=K\Norm_A(\Aut_T(R)).\]

\medskip

\noindent \textbf{Step 1}: \textit{We have $R\cap T=\Norm^K_Q(R)=\Norm_Q(R)$, and $R\Norm_Q(R)>R$.} The first equality is obvious. Let $g\in \Norm_Q(R)$, and $x\in R$. Then it is easy to see that $Qx^g=Qx$, so that $g$ acts trivially on $QR/Q$. Hence the automorphism determined by $g$ is in $K$, and so our first claim is proved. To see the second part, notice that $Q\cap R<Q$, and so $\Norm_Q(R)=\Norm_Q(Q\cap R)>Q\cap R$.

\medskip

\noindent \textbf{Step 2}: \textit{If $\psi\in \Norm_A(\Aut_T(R))$, then $N_\psi$ contains $T$.} Since $N_\psi$ is the inverse image under $\delta_R$ of the subgroup $\Aut_P(R)\cap \Aut_P(R)^{\psi^{-1}}$ in $A$, we need to show that $\Aut_T(R)$ is contained in both terms of the intersection. That it is contained in the first is clear, and for the second, since $\psi\in \Norm_A(\Aut_T(R))$, we have that $(\Aut_T(R))^\psi=\Aut_T(R)$. Thus it is contained in both terms, and so our claim is proved.

\medskip

Now we may prove the result: since $A=K\Norm_A(\Aut_T(R))$, the morphism $\phi$ may be written as $\phi=\chi\psi$, where $\chi\in K$ and $\psi\in \Norm_A(\Aut_T(R))$. Since $\chi$ acts trivially on $QR/Q$, we see that $\bar\phi=\bar\psi$ in $\bar\F_Q$. Furthermore, by Step 2 we see that $N_\psi$ contains $T$. However, if $N_\psi>R$, then $\psi$ extends to $\psi'$, on an overgroup of $R$, and so by induction $\bar\psi'$ lies in $\F/Q$. Since $\F/Q$ is a fusion system, this would imply that $\bar\psi$ is in $\F/Q$. Therefore $T\leq R$, and hence $Q\cap T=\Norm_Q(R)\leq R$. However, by Step 1, $R\Norm_Q(R)>R$, a contradiction, proving the theorem.
%
%\noindent \textbf{Step 3}: \textit{If $\psi\in \Norm_A(\Aut_T(R))$, then $\psi$ extends to an automorphism of $T$ that stabilizes $R$.} Here a miracle happens. (I see that $\psi$ extends to a map starting at $T$, and it stabilizes $R$.
%
%\medskip
%
%Since $\psi$ determines an element $\theta$ of $\Aut_\F(T)$, which normalizes $T\cap Q=\Norm_Q(R)$ clearly, we see that $\psi\in \Norm_A(\Aut_T(R))$ normalizes $\Aut_Q(R)$. Therefore
%\[\Norm_A(\Aut_T(R))\leq \Norm_A(\Aut_Q(R)),\]
%and since $\Inn(R)=\Aut_R(R)$ is a normal subgroup of $A$, we get that
%\[ \Norm_A(\Aut_T(R))\leq\Norm_A(\Aut_S(R))\]
%where $S=R(T\cap Q)>R$. Going back to our earlier Frattini argument, we see now that
%\[ A=K\Norm_A(\Aut_T(R))=K\Norm_A(\Aut_S(R)).\]
%
%\medskip
%
%\noindent \textbf{Step 4}: \textit{There are elements $\chi\in K$ and $\psi\in \Hom_\F(R,P)$ such that $\psi$ extends to $S$ and $\phi=\chi\psi$.} In the same way that we proved that $T\leq N_\psi$ for any $\psi\in \Norm_A(\Aut_T(R))$, we see that $S\leq N_\psi$ for any $\psi\in \Norm_A(\Aut_S(R))$. Since $\phi\in \Aut_\F(R)$, we see that there exists $\chi\in K$ and $\psi\in\Norm_A(\Aut_S(R))$ such that $\phi=\chi\psi$.
%
%\medskip
%
%Since $\chi$ acts trivially on $QR/Q$, we have that $\phi$ and $\psi$ have the same image in $\bar\F_Q$, and $\psi$ extends to a morphism in $\Hom_\F(S,P)$. Since $S>R$, by induction the image of $\psi$ lies in $\F/Q$, and so the image of $\phi$ lies in $\F/Q$, as needed.
\end{pf}

Thus if $\F$ is a saturated fusion system and $Q$ is a strongly $\F$-closed subgroup, then one may use either of the systems $\F/Q$ or $\bar\F_Q$ when making arguments about quotient systems. This will be essential in our short proof of Theorem \ref{productstronglyclosed}.

At this point we need to make a remark about the factor system $\F/Q$; Markus Linckelmann, in private communication, has pointed out that $\F/Q=\Norm_\F(Q)/Q$, and so $\F/Q$ is determined \emph{locally}. The interaction between this statement and Theorem \ref{saturatedquotient} might have significant implications for the structure of fusion systems, which have not yet been considered.

We end this section with what are essentially the second and third isomorphism theorems for fusion systems. Recall the definition of the universal fusion system on a finite $p$-group, consisting of all injective homomorphisms between subgroups of $P$.

\begin{prop}[Second Isomorphism Theorem]\label{secondisothm} Let $\F$ be a saturated fusion system on a finite $p$-group $P$, and let $Q$ be a strongly $\F$-closed subgroup, and let $\E$ be a saturated subsystem on a subgroup $R$. Write $\E Q/Q$ for the image of $\E$ in $\bar\F_Q$. Then
\[ \E Q/Q\cong \E/(R\cap Q).\]
\end{prop}
\begin{pf} The isomorphism $RQ/Q\to R/R\cap Q$ induces an isomorphism $\Phi:\mc U(RQ/Q)\to \mc U(R/R\cap Q)$ of the universal fusion systems, so we need to prove that the image of $\E Q/Q$ in $\mc U(R/R\cap Q)$ lies inside $\E/R\cap Q$ and vice versa.

Let $S/Q$ and $T/Q$ be subgroups of $RQ/Q$. Then a morphism $\phi:S/Q\to T/Q$ lies in $\E Q/Q$ if and only if there exist subgroups $S'$ and $T'$ of $R$ with $S'Q=S$ and $T'Q=T$, and a morphism $\psi\in\Hom_\E(S',T')$ such that the image of $\psi$ in $\bar\F_Q$ is $\phi$. The image of $\psi$ in $\E/R\cap Q$ is clearly $\phi\Phi$, and so the image of $\E Q/Q$ under $\Phi$ is contained in $\E/R\cap Q$. Conversely, if $\theta:S/R\cap Q\to T/R\cap Q$ is a morphism in $\E/R\cap Q$, then there is a morphism $\chi:S\to T$ with image $\theta$ in $\E/R\cap Q$, and the image $\bar\chi$ of $\chi$ in $\E Q/Q$ also satisfies $\bar\chi\Phi=\theta$, and so $\Phi$ induces an isomorphism $\E Q/Q\to \E/R\cap Q$, as needed.
\end{pf}

We get the following corollary \emph{a posteriori}.

\begin{cor} Let $\F$ be a saturated fusion system, and let $Q$ be a strongly $\F$-closed subgroup. The image of any saturated subsystem of $\F$ in $\bar\F_Q$ is saturated.
\end{cor}

We now consider the third isomorphism theorem for fusion systems.

\begin{prop}[Third Isomorphism Theorem]\label{thirdisothm} Let $\F$ be a saturated fusion system on a finite $p$-group $P$, and let $Q$ and $R$ be strongly $\F$-closed subgroups with $Q\leq R$. Then
\[ (\F/Q)/(R/Q)\cong \F/R.\]
\end{prop}
\begin{pf} By the third isomorphism theorem for groups, the two fusion systems $\E=(\F/Q)/(R/Q)$ and $\F/R$ are on the same subgroup. Suppose that $\bar\phi:S/R\to T/R$ is a morphism in $\F/R$. Then there is some morphism $\phi\in\Hom_\F(S,T)$ with image $\bar\phi$. Furthermore, the image $\phi':S/Q\to T/Q$ of $\phi$ in $\F/Q$ has image $\phi'':S/R\to T/R$, and since both $\bar\phi$ and $\phi''$ are derived from $\phi$, they must be the same morphism. The converse is a similar calculation, and is safely omitted.
\end{pf}

It is not clear whether the subsystem generated by two (particularly normal) subsystems is the right object to consider the product of the subsystems. We would like the product of two normal subsystems to itself be normal; since the product of two strongly closed subgroups is strongly closed, there is a possibility of this being true.

It is also not clear whether the subsystem $\E Q$ -- the full preimage of $\E Q/Q$ in $\F$ -- is saturated. If $\mc D$ is a normal subsystem on $Q$, then we would want $\E\mc D$ to be contained within $\E Q$ and contain both $\E$ and $\mc D$; any subsystem with these properties would satisfy the second isomorphism theorem, by Proposition \ref{secondisothm}.

\section{Closure and Quotients}

We begin with how weak and strong closure relates to taking quotients. The proof of (iv) in \cite{stancu2006} is not clear, because it is only true given Theorem \ref{saturatedquotient}, a result that is not mentioned in \cite{stancu2006}.

\begin{thm}[Stancu {{\cite[Lemma 4.7]{stancu2006}}}, Aschbacher {{\cite[Lemma 8.9]{aschbacher2008}}}]\label{weakstrongclosurefactor} Let $\F$ be a fusion system on a finite $p$-group $P$, and let $Q$ be a strongly $\F$-closed subgroup of $P$, and let $R$ be a subgroup of $P$.
\begin{enumerate}
\item The map $\Phi:\F\to \F/Q$ induces a bijection between the weakly $\F$-closed subgroups of $P$ containing $Q$ and those of $\F/Q$.
\item If $R$ is weakly $\F$-closed then the image of $R$ in $\F/Q$ is weakly $\F$-closed.
\item If $\F$ is saturated then the map $\Phi:\F\to \F/Q$ induces a bijection between the strongly $\F$-closed subgroups of $P$ containing $Q$ and those of $\F/Q$.
\item If $\F$ is saturated and $R$ is strongly $\F$-closed then the image of $R$ in $\F/Q$ is strongly $\F$-closed.
\end{enumerate}
\end{thm}
\begin{pf} Let $R$ be a subgroup containing $Q$. Then $R$ is weakly closed if and only if any morphism $\phi:R\to P$ is an automorphism. Clearly if $R$ is weakly $\F$-closed then any morphism $\phi':R/Q\to P/Q$ in $\F/Q$ (which must have a preimage in $\Hom_\F(R,P)$) is an automorphism, and vice versa. Hence (i) is true.

Also, the product of two weakly $\F$-closed subgroups is weakly $\F$-closed, so if any subgroup $R$ is weakly $\F$-closed, so is $QR$, and hence $QR/Q$ is weakly $\F/Q$-closed by (i), proving (ii).

For the rest of the proof, suppose that $\F$ is saturated, so that $\bar\F_Q=\F/Q$. Let $R$ be an overgroup of $Q$, and let $S$ be any subgroup of $R$. Suppose that $\phi:S\to P$ is a morphism in $\F$. Then $S\phi\leq R$ if and only if $(SQ/Q)\bar\phi=(S\phi)Q/Q\leq R/Q$, where $\bar\phi$ is the image in $\bar\F_Q$ of $\phi$. Since this holds for all $\phi$, we see that $R$ is strongly $\F$-closed if and only if $R/Q$ is strongly $(\bar\F_Q=\F/Q)$-closed, proving (iii).

Let $R$ be any strongly $\F$-closed subgroup of $P$, let $S$ be any subgroup of $R$, and let $\phi\in\Hom_\F(SQ,P)$. Since $Q$ is strongly $\F$-closed, $\phi|_Q$ is an automorphism of $Q$, and since $R$ is strongly $\F$-closed, $\phi|_S$ maps $S$ to $R$. Therefore, $(SQ)\phi\leq RQ$. If every subgroup of $RQ/Q$ is of the form $SQ/Q$ for some $S\leq R$, then $RQ/Q$ would be strongly $\F/Q$-closed, proving (iv). However, this follows from the second isomorphism theorem, since $RQ/Q\cong R/R\cap Q$, choosing $A/Q$ corresponds to a subgroup $B/R\cap Q$ on the right-hand side, and $BQ/Q=A/Q$, as needed.\end{pf}

We now prove Theorem \ref{normaliffstrongcentral}, by collating the equivalent conditions to a subgroup $Q$ having the property that $\F=\Norm_\F(Q)$.

\begin{thm} Let $\F$ be a saturated fusion system on a finite $p$-group $P$, and let $Q$ be a subgroup of $P$. The following are equivalent:
\begin{enumerate}
\item $\F=\Norm_\F(Q)$;
\item $\F_Q(Q)\normal \F$;
\item $Q$ is contained in every fully normalized, $\F$-centric, $\F$-radical subgroup of $\F$;
\item there is a central series for $Q$ all of whose terms (including $Q$) are strongly $\F$-closed; and
\item there is a central series for $Q$ all of whose terms are weakly $\F$-closed, and $Q$ is strongly $\F$-closed.
\end{enumerate}
\end{thm}
\begin{pf} The equivalence of (i) and (ii) is Theorem \ref{normaliffnormal}, the equivalence of (i) and (iii) is Proposition \ref{normalifffrc}, the equivalence of (i) and (iv) is Proposition \ref{normaliffstrongclosedcentral}, and that (iv) implies (v) is obvious. It remains to show that (v) implies (iv).

Let
\[ 1=Q_0\leq Q_1\leq \cdots \leq Q_n=Q\]
be a central series for $Q$, all of whose terms are weakly $\F$-closed. By Lemma \ref{weakcentralimpstrong}, $Q_1\leq \centre Q$ is strongly $\F$-closed. Since $Q_1$ is strongly $\F$-closed, we may take the quotient system $\F/Q_1$. Theorem \ref{weakstrongclosurefactor} states that the map $\F\to\F/Q_1$ induces a bijection between the weakly and strongly $\F$-closed subgroups of $P$ containing $Q_1$ and those of $P/Q_1$.

At this stage, one may either proceed by induction on the length of a central series all of whose terms are weakly $\F$-closed, by noticing that now $Q_2/Q_1$ is strongly $\F/Q_1$-closed, and hence $Q_2$ is strongly $\F$-closed, or proceed by induction on $|Q|$, and note that since $Q/Q_1$ has a satisfies (iv) of the theorem, there must be a central series by (v) which, when full preimages are taken, gives a central series all of whose terms are strongly $\F$-closed. Hence all conditions are equivalent, as claimed.
\end{pf}

Using Theorem \ref{weakstrongclosurefactor}, the proof of Theorem \ref{productstronglyclosed} is trivial. Let $Q$ and $R$ be strongly $\F$-closed subgroups of a saturated fusion system $\F$. By Theorem \ref{weakstrongclosurefactor}(iv), $QR/Q$ is strongly $\F/Q$-closed, and by (iii) of that theorem, this implies that $QR$ is strongly $\F$-closed, as required.

\bigskip

We now include a couple of lemmas, ready for our treatment of $p$-soluble fusion systems.

\begin{lem}\label{centrelift} Let $\F$ be a saturated fusion system on a finite $p$-group $P$, and suppose that $Z$ is a strongly $\F$-closed subgroup of $\centre\F$. Then for any normal subgroup $Q$, $\F=\Norm_\F(Q)$ if and only if $\F/Z=\Norm_{\F/Z}(Q/Z)$.
\end{lem}
\begin{pf} This is a special case of \cite[Theorem 6.5]{linckelmann2007}.\end{pf}

\begin{lem}\label{centnormalnorm} Let $\F$ be a saturated fusion system on a finite $p$-group $P$, and let $Q$ be a fully normalized subgroup. Let $K$ be a normal subgroup of $\Aut_\F(Q)$; then $\Norm_\F^K(Q)\normal \Norm_\F(Q)$. In particular, $\Cent_\F(Q)\normal \Norm_\F(Q)$.
\end{lem}
\begin{pf} Since $\F$ is saturated and $Q$ is fully normalized, $\Norm_\F(Q)$ is saturated. Also, since $K$ is a normal subgroup of $\Aut_\F(Q)$, $\Norm_\F^K(Q)$ is saturated as well, as we remarked when we defined $\Norm_\F^K(Q)$.

Next, we need to show that $\Norm_P^K(Q)$ is strongly $\Norm_\F(Q)$-closed, so let $R$ be a subgroup of $N=\Norm_P^K(Q)$ and $\phi:R\to S$ be a morphism in $\Norm_\F(Q)$. For $g\in R$, consider the action of $\theta_g$ on $Q$; since $R\leq N$, we see that $\theta_g\in K$. Furthermore, the action of $g\phi\in S$ on $Q$ is given by 
\[\theta_{g\phi}=\phi^{-1}\theta_g\phi\in K,\]
since $K\normal\Aut_\F(Q)$. Therefore $S\leq\Norm_P^K(Q)$, and so $N$ is strongly $\F$-closed.

Finally, we need to show that $\Norm_\F^K(Q)$ is $\Norm_\F(Q)$-invariant, so let $R\leq S\leq N$, $\phi\in\Hom_{\Norm_\F^K(Q)}(R,S)$ and $\psi\in\Hom_{\Norm_\F(Q)}(S,N)$. Since each of $\phi$ and $\psi$ extends to maps $\bar\phi$ and $\bar\psi$ whose domains include $Q$, and in the first case $\bar\phi|_Q\in K$ and in the second $\bar\phi|_Q\in\Aut_\F(Q)$, we see that
\[\(\bar\psi^{-1}\bar\phi\bar\psi\)|_Q\in K,\]
as $K$ is a normal subgroup of $\Aut_\F(Q)$. Hence $\psi^{-1}\phi\psi$ extends to a map $\theta$ whose domain includes $Q$ and for which $\theta|_Q\in K$. Therefore $\psi^{-1}\phi\psi\in\Norm_\F^K(Q)$, and so $\Norm_\F^K(Q)$ is $\Norm_\F(Q)$-invariant, as required.
\end{pf}

\section{$p$-Soluble Fusion Systems}
\label{psoluble}

For finite groups, a group is called $p$-soluble if repeated quotienting by (alternating) $\Orth_p(G)$ and $\Orth_{p'}(G)$ reaches the identity. In the case of fusion systems, we have that $\Orth_{p'}(\F)=1$, so it makes sense to make the following definition.

\begin{defn} Let $\F$ be a saturated fusion system on a finite $p$-group $P$. We say that $\F$ is \emph{$p$-soluble} if there exists a chain of strongly $\F$-closed subgroups
\[ 1=P_0\leq P_1\leq \cdots \leq P_n=P,\]
such that $P_i/P_{i-1}\leq \Orth_p(\F/P_{i-1})$ for all $1\leq i\leq n$. If $\F$ is $p$-soluble, then the length $n$ of a smallest such chain above will be called the \emph{$p$-length} of $\F$.
\end{defn}

The following lemma describes the basic facts of $p$-soluble fusion systems, mirroring those of finite groups. Define $\Orth_p^{(0)}(\F)=1$, and the $i$th term by
\[ \Orth_p^{(i)}(\F)/\Orth_p^{(i-1)}(\F)=\Orth_p\(\F/\Orth_p^{(i-1)}(\F)\).\]

\begin{lem}\label{psolublesubsystems} Let $\F$ be a saturated fusion system on a finite $p$-group $P$, let $Q$ be a strongly $\F$-closed subgroup of $P$ and let $\E$ be a normal subsystem of $\F$.
\begin{enumerate}
\item If $\F$ is $p$-soluble then all saturated subsystems and quotients $\F/Q$ are $p$-soluble.
\item If $\E$ and $\F/\E$ are $p$-soluble then so is $\F$.
\item $\F$ is $p$-soluble if and only $\Orth_p^{(n)}(\F)=P$ for some $n$, and the smallest such $n$ is the $p$-length of $\F$.
\end{enumerate}
\end{lem}
\begin{pf} Choose a triple $(\F,P,Q)$, where $\F$ is a $p$-soluble, saturated fusion system on a finite $p$-group $P$, and $Q$ is a strongly $\F$-closed subgroup of $P$ such that $\F/Q$ is not $p$-soluble, and such that $|P|$ is minimal subject to these constraints. Write $R=\Orth_p(\F)\neq 1$, and we claim that $QR/Q\leq \Orth_p(\F/Q)$. To see this notice that any morphism $\phi:S\to T$ extends to a morphism $\psi:RS\to RT$, and the image $\bar\phi$ of $\phi$ in $\bar\F_Q$ is extended by the image $\bar\psi$ of the morphism $\psi$ in $\bar\F_Q$. Since $\F/Q$ is not $p$-soluble, neither is (using the third isomorphism theorem)
\[ \left.(\F/Q)\right/(QR/Q)\cong \F/QR.\]
However, clearly $\F/\Orth_p(\F)$ is $p$-soluble, and so $(\F/R,P/R,QR/R)$ is a triple satisfying our conditions with $|P|>|P/R|$. This yields a contradiction, proving that quotients of $p$-soluble fusion systems are $p$-soluble.

Now let $(\F,\F',Q)$ be a triple, with $\F$ a $p$-soluble fusion system on a $p$-group $P$, $\F'$ a saturated subsystem of $\F$ on a subgroup $Q$, with $\F'$ not $p$-soluble. Choose this triple with $|P|$ minimal. Let $R=\Orth_p(\F)\neq 1$, and we claim that $R\cap Q\leq \Orth_p(\F')$. Since $\Norm_\F(R)=\F$, we have a central series
\[ 1=R_0\leq R_1\leq \cdots \leq R_n=R,\]
with each $R_i$ strongly $\F$-closed, by Proposition \ref{normaliffstrongclosedcentral}. Therefore any morphism in $\E$ whose domain lies inside $R_i$ has image inside $R_i$. Since $Q$ is obviously strongly $\E$-closed, the means that the intersection $Q_i=R_i\cap Q$ is strongly $\E$-closed. Thus the series of the $Q_i$ is a central series for $Q\cap R$ whose terms are strongly $\E$-closed, which by another application of Proposition \ref{normaliffstrongclosedcentral}, gives the result.

Consider the image $\F' R/R$ of $\F'$ in $\F/R$. Since this is isomorphic with $\F'/R\cap Q$ by the second isomorphism theorem, and is hence not $p$-soluble, the triple $(\F/R,\F' R/R,QR/R)$ has $|P/R|<|P|$ and so this contradicts the choice of the original triple. This proves (i).

If $\E$ is $p$-soluble, then $\mc O_p(\E)\charac \E\normal \F$, so that $\Orth_p(\E)$ is strongly $\F$-closed; this, together with induction and Theorem \ref{weakstrongclosurefactor}, proves that $\Orth_p^{(i)}(\E)$ is strongly $\F$-closed for all $i$. This includes $R$, one of the $\Orth_p^{(i)}(\E)$, and so by Theorem \ref{weakstrongclosurefactor} again, the preimages of $\Orth_p^{(j)}(\F/R)$ are strongly $\F$-closed. The concatenation of the two series $\Orth_p^{(i)}(\E)$ and the preimages of $\Orth_p^{(j)}(\F/R)$ satisfy the requirement for $\F$ to be $p$-soluble, proving (ii).

Now, we prove (iii). The one direction is clear, so suppose that $\F$ has a series $(Q_i)$ so that $Q_i/Q_{i-1}\leq \Orth_p(\F/Q_{i-1})$ for all $i$. If we can show that, whenever $S\leq \Orth_p(\F)$ is strongly $\F$-closed, then the preimage of $\Orth_p(\F/S)$ is contained in $\Orth_p^{(2)}(\F)$, then we are done by an obvious induction.

We notice that if $S\leq \Orth_p(\F)=R$, then $\Orth_p(\F/S)$ contains $R$; to see this, a morphism $\bar\phi:A/S\to B/S$ in $\F/S$ comes from a morphism $\phi:A\to B$ in $\F$, and this extends to a morphism $\psi:AR\to BR$, which has image $\bar\psi:AR/S\to BR/S$, extending $\bar\phi$. Therefore $R/S\leq \Orth_p(\F/S)$, and so one may consider $\Orth_p(\F/S)/R$, which is contained in $\Orth_p^{(2)}(\F)/R$ by another application of the claim. Thus all parts of the lemma are proved.
\end{pf}

At this point we digress briefly to discuss minimal normal subsystems. At the moment, there is no characterization of minimal normal subsystems, like there is for groups. The reason behind this is that the intersection of two normal subsystems need not be normal, nor even saturated.

\begin{example} Let $P=D_8\times C_2$, with the $D_8$ factor generated by an element $x$ of order $4$ and $y$ of order $2$, and the $C_2$-factor being generated by $z$. Let $Q=\gen{x,y}$, and $R=\gen{xz,y}$. Then $S=Q\cap R$ is a normal Klein four subgroup of $P$, and $\Aut_Q(S)=\Aut_R(S)$ contains the identity and the map swapping $y$ and $x^2y$. Thus $\E=\F_Q(Q)\cap \F_R(R)$ has an outer automorphism of order $2$ on $S$, and so cannot be saturated, as $\Aut_S(S)$ is not a Sylow $2$-subgroup of $\Aut_\E(S)$.
\end{example}

This example shows that for any `reasonable' definition of normality -- i.e., one for which normal subgroups yield normal subsystems -- the intersection of two normal subsystems need not be saturated. In \cite{aschbacher2007un}, Aschbacher proves  -- using the definition of normality given in \cite{aschbacher2008}, which we will call \emph{strong normality}, but not define here -- the intersection of two strongly normal subsystems \emph{contains} a strongly normal subsystem on the intersection of the two relevant subgroups. Such a result is not known for normal subsystems, but if it were true, then one could prove that any minimal normal subsystem is either $\F_Q(Q)$ for some elementary abelian $p$-group $Q$, or the direct product of isomorphic simple fusion systems. By mimicking the proof for groups, one can get this result for minimal, \emph{strongly} normal subsystems, but it needs the theorem on intersections mentioned above.

However, we can say something about minimal \emph{subnormal} subsystems. (The definition of subnormality is obvious, and left to the reader.) These are obviously simple, and so either of the form $\F_Q(Q)$ for $Q$ of prime order, or some non-abelian simple fusion system. The point here is the following.

\begin{prop}\label{solublesimple} Let $\E$ be a minimal subnormal subsystem of the saturated fusion system $\F$, on a finite $p$-group $P$. If $\E=\F_Q(Q)$ for $Q\leq P$ of prime order, then $Q\leq \Orth_p(\F)$. In particular, $\F$ is not $p$-soluble if and only if there is a strongly $\F$-closed subgroup $R$ and non-abelian simple subsystem of $\F/R$.\end{prop}
\begin{pf} The first statement, that $\Orth_p(\E)\leq \Orth_p(\F)$ for any subnormal subsystem $\E$ of $\F$, easily follows from Proposition \ref{normalop} and induction on the length of a chain of normal subsystems connecting $\E$ and $\F$.

To see the second, if $\F$ is $p$-soluble then all quotients and saturated subsystems are $p$-soluble, and so $\F$ contains no non-abelian simple subquotient, proving one direction. For the other, proceed by induction on $|P|$; assume that $\F$ is a saturated fusion system on a finite $p$-group $P$ with no non-abelian simple subquotients. If $\Orth_p(\F)\neq 1$, then $\F/\Orth_p(\F)$ is $p$-soluble by induction, whence so is $\F$. However, no minimal subnormal system may be non-abelian simple, and so must be of the form $\F_Q(Q)$, and so we see that $\Orth_p(\F)\neq 1$ by the first part of the proposition.
\end{pf}

We come to the main result of the section, the proof that $p$-soluble fusion systems are constrained. We start with a simple lemma.

\begin{lem}\label{opequalz} For any saturated fusion system $\F$, $\centre\F\leq \Orth_p(\F)$, and
\[ \Orth_p(\F)/\centre\F=\Orth_p(\F/\centre\F).\]
In particular, if $\Orth_p(\F)=\centre\F$ then either $P=\Orth_p(\F)$ or $\F$ is not $p$-soluble.
\end{lem}
\begin{pf} Clearly $\centre\F\leq \Orth_p(\F)$, and by Lemma \ref{centrelift}, 
\[ \Orth_p(\F)/\centre\F=\Orth_p(\F/\centre\F).\]
If $\Orth_p(\F)\neq P$, then $\Orth_p(\F/\centre\F)$ is trivial. Therefore no minimal subnormal subsystem $\E$ of $\F/Q$ is of the form $\F_Q(Q)$ for some $Q$, and so Proposition \ref{solublesimple} implies that $\F$ is not $p$-soluble, as claimed.
\end{pf}

\begin{thm} Suppose that $\F$ is $p$-soluble, and let $Q=\Orth_p(\F)$. Then $\Cent_P(Q)=\centre Q$, so that $\F$ is constrained.
\end{thm}
\begin{pf} Let $Q=\Orth_p(\F)$, and let $\E=\Cent_\F(Q)$. Since $\F=\Norm_\F(Q)$, we see that $\E\normal\F$ by Lemma \ref{centnormalnorm}. Therefore, by Proposition \ref{normalop}, 
\[ \Orth_p(\E)=\Orth_p(\F)\cap \Cent_P(Q)=Q\cap \Cent_P(Q)=\centre Q.\]

However, since $\E=\Cent_\F(Q)$, we notice that every morphism in $\E$ centralizes $\centre Q=Q\cap\Cent_P(Q)$, so that $\centre Q\leq \centre \E$. Since $\centre \E\leq \Orth_p(\E)$ obviously, we see that
\[ \Orth_p(\E)=\centre \E.\]
As $\E$ is $p$-soluble (since $\F$ is, by Lemma \ref{psolublesubsystems}), we get by Lemma \ref{opequalz} that $\Cent_P(Q)=\Orth_p(\Cent_\F(Q))$. However, $\Orth_p(\E)\leq \Orth_p(\F)$ by Proposition \ref{normalop}, so that
\[ \Cent_P(Q)=\Orth_p(\E)\leq \Orth_p(\F)=Q.\]
Thus $Q$ is $\F$-centric, as required.
\end{pf}

\section{Soluble Systems and Soluble Groups}

The fact that $p$-soluble fusion systems are constrained means that they are fusion systems of finite groups, by Theorem \ref{constrainedmodels}. In this section we will characterize those fusion systems that come from \emph{$p$-soluble} groups. Since every $p$-soluble group has a $p$-soluble fusion system, it is a necessary condition that they be constrained, so we start from here.

It is well-known (see, for example, \cite[33.12]{aschbacher}) that if $G$ is a simple $p'$-group, then any quasisimple group with quotient $G$ is also a $p'$-group. This gives us the following, originally proved by Hall and Higman.

\begin{lem}\label{psolimpconst} Let $G$ be a $p$-soluble group with $\Orth_{p'}(G)=1$. Then
\[ \Cent_G(\Orth_p(G))\leq \Orth_p(G).\]
\end{lem}
\begin{pf} A fundamental property of the generalized Fitting subgroup is that $\Cent_G(F^*(G))\leq F^*(G)$ \cite[31.13]{aschbacher}. Since $G$ is $p$-soluble, any subnormal quasisimple group is a $p'$-group, and therefore $E(G)\leq \Orth_{p'}(G)=1$. Therefore $F^*(G)=\Orth_p(G)$, proving the claim.\end{pf}

Using this lemma, we give our first characterization of fusion systems of $p$-soluble groups. In the rest of this section, for a constrained fusion system $\F$, we denote the unique finite group given in Theorem \ref{constrainedmodels} by $L_\F$.

\begin{prop}\label{soliffconstsol} Let $\F$ be a saturated fusion system on a finite $p$-group $P$. There is a $p$-soluble group $G$ such that $\F=\F_P(G)$ if and only if $\F$ is constrained, and $L_\F$ is $p$-soluble.
\end{prop}
\begin{pf} Suppose that $\F=\F_P(G)$ for some $p$-soluble group $G$. We may assume that $\Orth_{p'}(G)=1$, since the fusion systems on $G$ and $G/\Orth_{p'}(G)$ are the same. By Lemma \ref{psolimpconst}, $\Orth_p(G)$ contains its centralizer. Therefore, by the uniqueness of the group in Theorem \ref{constrainedmodels}, $G=L_\F$.
\end{pf}

Finally, this allows us to reach an internal characterization of fusion systems of $p$-soluble groups, without reference to groups.

\begin{cor}Let $\F$ be a saturated fusion system on a finite $p$-group $P$, and write $Q=\Orth_p(\F)$. Then $\F$ is the fusion system of a $p$-soluble group if and only if $\F$ is constrained and $\Aut_\F(Q)$ is $p$-soluble.
\end{cor}
\begin{pf} Suppose that $\F$ is constrained and that $\Aut_\F(Q)$ is $p$-soluble. Since $\F$ is constrained, we have that $\F=\F_P(L_\F)$; then $L_\F$ is an extension of $\Cent_{L_\F}(Q)=\centre Q$ by $\Aut_{L_\F}(Q)=\Aut_\F(Q)$, which is $p$-soluble. Hence $L_\F$ is $p$-soluble, and so $\F$ is the fusion system of a $p$-soluble group. The converse is similarly clear.
\end{pf}

At first blush, this appears to contradict an assertion in \cite[1.6]{puig2006}, which claims that all $p$-soluble fusion systems arise from $p$-soluble groups. The incongruity stems from the definition of a soluble fusion system in \cite{puig2006}, which works `from the top down', in the sense that it involves taking repeated subsystems, rather than repeated quotients. Recall the definitions of $\Orth^p(\F)$ and $\Orth^{p'}(\F)$ from \cite{puig2006}, which will not be repeated here. Then another way to define a soluble fusion system is that repeated taking of $\Orth^p$ and $\Orth^{p'}$ operators eventually reaches the trivial group. While this definition picks out exactly the fusion systems of $p$-soluble groups (this is not difficult to prove given Theorem \ref{intcharsol}), it suffers from the fact that this class of fusion system is not closed under extensions. The class of $p$-soluble fusion systems given here is extension-closed and contains the fusion system $\F_P(P)$, where $P$ is cyclic of order $p$; by Proposition \ref{solublesimple}, it is also the class of all systems that do not contain simple subquotients. It seems to us that this is the `correct' definition of solubility for fusion systems, given the use of the word `soluble' throughout algebra to mean similar concepts.

\bigskip

We now consider so-called $Qd(p)$-free fusion systems. Firstly, the group $Qd(p)$ is the semidirect product $(C_p\times C_p)\rtimes \SL_2(p)$, with $\SL_2(p)$ acting in the natural way. A finite group is called \emph{$Qd(p)$-free} if no subquotient of it is isomorphic with $Qd(p)$. As in \cite{kessarlinckelmann2008}, define a saturated fusion system $\F$ to be \emph{$Qd(p)$-free} if, for any fully normalized, $\F$-centric subgroup $Q$ of $P$, The group $L_{\Norm_\F(Q)}$ is $Qd(p)$-free.

\begin{lem} Every $Qd(p)$-free fusion system $\F$ is $p$-soluble.
\end{lem}
\begin{pf} By \cite[Proposition 6.4]{kessarlinckelmann2008}, if $\F$ is $Qd(p)$-free then so is $\F/\Orth_p(\F)$, so it suffices to show that $\Orth_p(\F)\neq 1$; this is given to us by 7.1 of the same paper.
\end{pf}

Since a $Qd(p)$-free fusion system is $p$-soluble it is constrained, and hence is the fusion system of some, $Qd(p)$-free, group. Thus any theorem known to hold for $Qd(p)$-free groups should have an analogue for $Qd(p)$-free fusion systems. One such example is Glauberman's ZJ-Theorem: another is Thompson Factorization.

\begin{thm}[Thompson Factorization] Let $\F$ be a $Qd(p)$-free fusion system on a finite $p$-group $P$, where $p$ is odd. Then
\[ \F=\Norm_\F(J(P))\Cent_\F(\Omega_1(\centre{P})).\]
\end{thm}
\begin{pf} Such a decomposition holds for $Qd(p)$-free groups (see, for example, \cite[Theorem 26.9]{glsvol2}). The rest is simply rewriting of the factorization into fusion systems.
\end{pf}

However, a na\"ive rewriting of Thompson Factorization for other cases fails; for example, Glauberman proved that if $p\geq 5$, then the conclusion to the above theorem holds for all $p$-soluble groups, with a slight modification to the centralizer term (see \cite[Theorem 26.10]{glsvol2}). This does not carry over to $p$-soluble fusion systems in general; for instance, the groups $Qd(p)$ themselves do not satisfy this version of Thompson Factorization.

\bigskip \bigskip

\noindent \textbf{Acknowledgments}: I would like to thank Adam Glesser for being a human sounding board for many if not all of the ideas in this paper, for pointing out some errors in the original version, and for suggesting improvements to the exposition. I also had a helpful conversation with Inna Korchigina on $2$-soluble fusion systems, in which we proved Proposition \ref{soliffconstsol} for this case.

\bibliography{references}

\end{document}